\documentclass[12pt,reqno]{article}

\usepackage[usenames]{color}
\usepackage{amssymb}
\usepackage{graphicx}
\usepackage{amscd}

\usepackage{amsthm}
\newtheorem{theorem}{Theorem}

\theoremstyle{definition}

\newtheorem{example}[theorem]{Example}

\usepackage[colorlinks=true,
linkcolor=webgreen, filecolor=webbrown,
citecolor=webgreen]{hyperref}

\definecolor{webgreen}{rgb}{0,.5,0}
\definecolor{webbrown}{rgb}{.6,0,0}

\usepackage{color}

\usepackage{float}

\usepackage{graphics,amsmath,amssymb}
\usepackage{amsfonts}
\usepackage{latexsym}
\usepackage{epsf}

\newcommand{\seqnum}[1]{\href{http://www.research.att.com/cgi-bin/access.cgi/as/~njas/sequences/eisA.cgi?Anum=#1}{\underline{#1}}}

\begin{document}

\begin{center}
\vskip 1cm{\LARGE\bf Comparing two matrices of generalized moments defined by continued fraction expansions} \vskip 1cm \large
Paul Barry\\
School of Science\\
Waterford Institute of Technology\\
Ireland\\
\href{mailto:pbarry@wit.ie}{\tt pbarry@wit.ie} \\

\end{center}
\vskip .2 in

\begin{abstract} We study two matrices $N$ and $M$ defined by the parameters of equivalent $S$- and $J$-continued fraction expansions, and compare them by examining the product $N^{-1}M$. Using examples based on the Catalan numbers, the little Schr\"oder numbers and powers of $q$, we indicate that this matrix product is an object worthy of study. In the case of the little Schr\"oder numbers, we find that the matrix $N$ has an interleaved structure based on two Riordan arrays.
\end{abstract}

\section{Introduction}

In this note, we study two matrices whose elements may be considered to be generalized moments. The matrices are defined using the coefficients of simple Jacobi and Stieltjes continued fractions.

In familiar cases, these matrices are well-known, though this study examines them from a fresh perspective. It will be assumed that the reader is familiar with the basics of orthogonal polynomials \cite{Chihara, Gautschi, Szego}, Riordan arrays \cite{SGWW}, production matrices \cite{Prod1, Prod2, PW}, continued fractions \cite{Wall} and the interplay between these areas \cite{Barry_Meixner, Barry_Moment}.

Our point of departure is a sequence $a_n$, with $a_1=1$, whose elements are either integers or polynomials with integer coefficients.

We will use these numbers to define two lower-triangular matrices, which we then compare.

In both cases, the elements of the first column will be the sequence $\mu_n$ generated by the continued fraction

$$\cfrac{1}{1-\cfrac{a_1 x}{1-\cfrac{a_2x}{1-\cfrac{a_3 x}{1-\cdots}}}}.$$

\noindent
We require that this sequence be Catalan-like, in the sense that we require all the Hankel determinants $|\mu_{i+j}|_{0 \le i,j \le n}$ to be non-zero.

An equivalence transformation ensures that the sequence $\mu_n$ is the same as that generated by
$$\cfrac{1}{1-a_1 x- \cfrac{a_1 a_2 x^2}{1-(a_2+a_3)x-\cfrac{a_3 a_4 x^2}{1-(a_4+a_5)x-\cdots}}}.$$
This exhibits $\mu_n$ as the moment sequence of the family of orthogonal polynomials $P_n(x)$ that satisfy
$$P_n(x)=(x-(a_{2n-2}+a_{2n-1}))P_{n-1}(x)-a_{2n-3}a_{2n-2}P_{n-2}(x),$$ with $P_0(x)=1$, $P_1(x)=x-a_1$. 

The first matrix $M$ that we shall be interested in is the inverse of the matrix of coefficients of these polynomials.
This matrix therefore has production matrix given by

\begin{displaymath}\left(\begin{array}{ccccccc} a_1 & 1 &
0
& 0 & 0 & 0 & \ldots \\a_1 a_2 & a_2+a_3 & 1 & 0 & 0 & 0 & \ldots \\ 0 & a_3 a_4
& a_4+a_5 & 1 & 0 &
0 & \ldots \\ 0 & 0 & a_5 a_6 & a_6+a_7 & 1 & 0 & \ldots \\ 0 & 0 & 0
& a_7 a_8 & a_8+ a_9 & 1 & \ldots \\0 & 0 & 0 & 0 & a_9 a_{10} & a_{10}+a_{11}
&\ldots\\
\vdots &
\vdots & \vdots & \vdots & \vdots & \vdots &
\ddots\end{array}\right).\end{displaymath}
\noindent The form of this production matrix ensures that the matrix generated by it will be lower-triangular with $1$'s on the diagonal. We obtain a matrix which begins

\begin{displaymath}\scriptsize\left(\begin{array}{cccccc}
1 & 0 & 0 & 0 & 0 & \ldots \\
\alpha & 1 & 0 & 0 & 0  & \ldots \\
\alpha(\alpha+\beta) & \alpha+\beta+\gamma & 1 & 0  & 0 & \ldots \\
\alpha((\alpha+\beta)^2+\beta \gamma) & (\alpha+\beta)^2+\beta \gamma+\gamma(\alpha+\beta+\gamma+\delta) & \alpha+\beta+\gamma+\delta+\epsilon  & 1 & 0 & \ldots \\

\alpha((\alpha+\beta)^3+\beta \gamma (2\alpha+2\beta+\gamma+\delta) & (\alpha+\beta+\gamma+\delta+\epsilon)^2-\alpha(\gamma+\delta+\epsilon)-\beta(\delta+\epsilon)-\epsilon(\gamma-\phi) & \cdots
& \alpha+\beta+\gamma+\delta+\epsilon+\phi+\psi & 1  & \ldots \\
\vdots &
\vdots & \vdots &  \vdots & \vdots &
\ddots\end{array}\right), \end{displaymath}
where we have written $\alpha=a_1$, $\beta=a_2$, and so on.

In order to define the second matrix $N$, which again will have $\mu_n$ in the first column, we also use a production matrix.
To construct this production matrix, we have two alternative routes. The first one proceeds as follows; we take the inverse of the matrix
\begin{displaymath}\left(\begin{array}{ccccccc} 1 & 0 &
0
& 0 & 0 & 0 & \ldots \\-a_1 & 1 & 0 & 0 & 0 & 0 & \ldots \\ 0 & -a_2
& 1 & 0 & 0 &
0 & \ldots \\ 0 & 0 & -a_3 & 1 & 0 & 0 & \ldots \\ 0 & 0 & 0
& -a_4 & 1 & 0 & \ldots \\0 & 0 & 0 & 0 & -a_5 & 1
&\ldots\\
\vdots &
\vdots & \vdots & \vdots & \vdots & \vdots &
\ddots\end{array}\right) \end{displaymath} to obtain the matrix

\begin{displaymath}\left(\begin{array}{ccccccc} 1 & 0 &
0
& 0 & 0 & 0 & \ldots \\a_1 & 1 & 0 & 0 & 0 & 0 & \ldots \\ a_1 a_2 & a_2
& 1 & 0 & 0 &
0 & \ldots \\ a_1 a_2 a_3 & a_2 a_3 & a_3 & 1 & 0 & 0 & \ldots \\ a_1 a_2 a_3 a_4 & a_2 a_3 a_4 & a_3 a_4
& a_4 & 1 & 0 & \ldots \\a_1 a_2 a_3 a_4 a_5 & a_2 a_3 a_4 a_5 & a_3 a_4 a_5 & a_4 a_5 & a_5 & 1
&\ldots\\
\vdots &
\vdots & \vdots & \vdots & \vdots & \vdots &
\ddots\end{array}\right). \end{displaymath}
We now behead this matrix (we remove the first row) to obtain the following production matrix.
\begin{displaymath}\left(\begin{array}{ccccccc} a_1 & 1 & 0 & 0 & 0 & 0 & \ldots \\ a_1 a_2 & a_2
& 1 & 0 & 0 &
0 & \ldots \\ a_1 a_2 a_3 & a_2 a_3 & a_3 & 1 & 0 & 0 & \ldots \\ a_1 a_2 a_3 a_4 & a_2 a_3 a_4 & a_3 a_4
& a_4 & 1 & 0 & \ldots \\a_1 a_2 a_3 a_4 a_5 & a_2 a_3 a_4 a_5 & a_3 a_4 a_5 & a_4 a_5 & a_5 & 1
&\ldots\\
a_1 a_2 a_3 a_4 a_5 a_6 & a_2 a_3 a_4 a_5 a_6 &
a_3 a_4 a_5 a_6
& a_4 a_5 a_6 & a_5 a_6 & a_6 & \ldots \\
\vdots &
\vdots & \vdots & \vdots & \vdots & \vdots &
\ddots\end{array}\right). \end{displaymath}
The form of this production matrix ensures that the matrix $N$ that it generates will be lower-triangular with $1$'s on the diagonal. The matrix $N$ that we seek then begins

\begin{displaymath}\scriptsize\left(\begin{array}{cccccc} 1 & 0 &
0
& 0 & 0 & \ldots \\ \alpha & 1 & 0 & 0 & 0  & \ldots \\ \alpha(\alpha+\beta) & \alpha+\beta
& 1 & 0  &
0 & \ldots \\ \alpha((\alpha+\beta)^2+\beta \gamma) & (\alpha+\beta)^2+\beta \gamma & \alpha+\beta+\gamma & 1 & 0 & \ldots \\

\alpha((\alpha+\beta)^3+\beta \gamma (2\alpha+2\beta+\gamma+\delta) & (\alpha+\beta)^3+\beta \gamma (2\alpha+2\beta+\gamma+\delta & (\alpha+\beta)^2+\gamma(\alpha+\delta)+(\beta+\gamma)^2
& \alpha+\beta+\gamma+\delta & 1  & \ldots \\
\vdots &
\vdots & \vdots &  \vdots & \vdots &
\ddots\end{array}\right), \end{displaymath} where we have used $\alpha=a_1$, $\beta=a_2$, and so on.

\noindent There is an alternative production matrix approach to the construction of $N$. Multiplying the $(n,k)$-th element of $N$ by
$$\prod_{j=1}^k a_j$$ produces a lower triangular matrix whose first column is the same as that of $N$, and whose production matrix takes the simple form of
\begin{displaymath}\left(\begin{array}{ccccccc} a_1 & a_1 &
0
& 0 & 0 & 0 & \ldots \\a_2 & a_2 & a_2 & 0 & 0 & 0 & \ldots \\ a_3 & a_3
& a_3 & a_3 & 0 &
0 & \ldots \\ a_4 & a_4 & a_4 & a_4 & a_4 & 0 & \ldots \\ a_5 & a_5 & a_5
& a_5 & a_5 & a_5 & \ldots \\a_6 & a_6 & a_6 & a_6 & a_6 & a_6
&\ldots\\
\vdots &
\vdots & \vdots & \vdots & \vdots & \vdots &
\ddots\end{array}\right).  \end{displaymath}
We can clearly reverse this process, starting with the sequence $a_n$, to produce $N$.

In order to compare the two matrices $M$ and $N$, it is natural to examine the product $N^{-1}M$.

\begin{example} \textbf{The Catalan matrices}. We let $a_n=1$. Thus we are interested in the sequence
generated by the continued fraction
$$\cfrac{1}{1-\cfrac{ x}{1-\cfrac{x}{1-\cfrac{ x}{1-\cdots}}}}.$$
\noindent This is the sequence of Catalan numbers $C_n = \frac{1}{n+1}\binom{2n}{n}$, \seqnum{A000108}. In this instance,
the production matrix of  $N$ for both methods of generation is given by
\begin{displaymath}\left(\begin{array}{ccccccc} 1 & 1 &
0
& 0 & 0 & 0 & \ldots \\1 & 1 & 1 & 0 & 0 & 0 & \ldots \\ 1 & 1
& 1 & 1 & 0 &
0 & \ldots \\ 1 & 1 & 1 & 1 & 1 & 0 & \ldots \\ 1 & 1 & 1
& 1 & 1 & 1 & \ldots \\1 & 1 & 1 & 1 & 1 & 1
&\ldots\\
\vdots &
\vdots & \vdots & \vdots & \vdots & \vdots &
\ddots\end{array}\right),  \end{displaymath} and $N$ is the Riordan array
$$N=(c(x), xc(x))=(1-x, x(1-x))^{-1} \quad\quad \seqnum{A033184}.$$
\noindent The matrix $M$ is given by the Riordan array
$$M=(c(x), xc(x)^2)=\left(\frac{1}{1+x}, \frac{x}{(1+x)^2}\right)^{-1} \quad\quad \seqnum{A039599},$$ and the associated orthogonal polynomials are the
Chebyshev polynomials $U_n(\frac{x}{2})$. The production matrix of $(c(x), xc(x)^2)$ is given by
\begin{displaymath}\left(\begin{array}{ccccccc} 1 & 1 &
0
& 0 & 0 & 0 & \ldots \\1 & 2 & 1 & 0 & 0 & 0 & \ldots \\ 0 & 1
& 2 & 1 & 0 &
0 & \ldots \\ 0 & 0 & 1 & 2 & 1 & 0 & \ldots \\ 0 & 0 & 0
& 1 & 2 & 1 & \ldots \\0 & 0 & 0 & 0 & 1 & 2
&\ldots\\
\vdots &
\vdots & \vdots & \vdots & \vdots & \vdots &
\ddots\end{array}\right),  \end{displaymath} corresponding to the generating function
$$\cfrac{1}{1-x-\cfrac{x^2}{1-2x-\cfrac{x^2}{1-2x-\cdots}}}$$ of $C_n$.

A straight-forward Riordan array calculation now shows that in this case,
$$N^{-1}\cdot M = (c(x), xc(x))^{-1} \cdot (c(x), xc(x)^2)= \left(1, \frac{x}{1-x}\right),$$ which is the shifted binomial matrix
\begin{displaymath}\left(\begin{array}{ccccccc} 1 & 0 &
0
& 0 & 0 & 0 & \ldots \\0 & 1 & 0 & 0 & 0 & 0 & \ldots \\ 0 & 1
& 1 & 0 & 0 &
0 & \ldots \\ 0 & 1 & 2 & 1 & 0 & 0 & \ldots \\ 0 & 1 & 3
& 3 & 1 & 0 & \ldots \\0 & 1 & 4 & 6 & 4 & 1
&\ldots\\
\vdots &
\vdots & \vdots & \vdots & \vdots & \vdots &
\ddots\end{array}\right).  \end{displaymath}

\end{example}
\begin{example} \textbf{The $q$-case}. We take the example of
$$a_n = \frac{q^n}{q}-\frac{0^n}{q}.$$
\noindent Starting with the matrix
\begin{displaymath}\left(\begin{array}{ccccccc} 1 & 0 &
0
& 0 & 0 & 0 & \ldots \\-q & 1 & 0 & 0 & 0 & 0 & \ldots \\ 0 & -q^2
& 1 & 0 & 0 &
0 & \ldots \\ 0 & 0 & -q^3 & 1 & 0 & 0 & \ldots \\ 0 & 0 & 0
& -q^4 & 1 & 0 & \ldots \\0 & 0 & 0 & 0 & -q^5 & 1
&\ldots\\
\vdots &
\vdots & \vdots & \vdots & \vdots & \vdots &
\ddots\end{array}\right), \end{displaymath} we invert it and behead the resulting matrix to get the production matrix
\begin{displaymath}\left(\begin{array}{ccccccc} 1 & 1 &
0
& 0 & 0 & 0 & \ldots \\q & q & 1 & 0 & 0 & 0 & \ldots \\ q^3 & q^3
& q^2 & 1 & 0 &
0 & \ldots \\ q^6 & q^6 & q^5 & q^3 & 1 & 0 & \ldots \\ q^{10} & q^{10} & q^9
& q^7 & q^4 & 1 & \ldots \\q^{15} & q^{15} & q^{14} & q^{12} & q^9 & q^5
&\ldots\\
\vdots &
\vdots & \vdots & \vdots & \vdots & \vdots &
\ddots\end{array}\right),  \end{displaymath} which we use to generate the matrix N:
\begin{displaymath}\scriptsize\left(\begin{array}{cccccc} 1 & 0 &
0
& 0 & 0 & \ldots \\1 & 1 & 0 & 0 & 0 & \ldots \\ q+1 & q+1
& 1 & 0 & 0 &
 \ldots \\ q^3+q^2+2q+1 & q^3+q^2+2q+1 & q^2+q+1 & 1 & 0 &  \ldots \\ q^6+q^5+2q^4+3q^3+3q^2+3q+1 & q^6+q^5+2q^4+3q^3+3q^2+3q+1 & q^5+q^4+2q^3+2q^2+2q+1
& q^3+q^2+q+1 & 1 &  \ldots \\
\vdots &
\vdots & \vdots & \vdots & \vdots &
\ddots\end{array}\right). \end{displaymath}
\noindent In the left column we recognize the $q$-Catalan sequence $\mu_n$ with generating function
$$\cfrac{1}{1-\cfrac{ x}{1-\cfrac{qx}{1-\cfrac{q^2 x}{1-\cdots}}}}.$$

\noindent Alternatively we may begin with the production matrix
\begin{displaymath}\left(\begin{array}{ccccccc} 1 & 1 &
0
& 0 & 0 & 0 & \ldots \\q & q & q & 0 & 0 & 0 & \ldots \\ q^2 & q^2
& q^2 & q^2 & 0 &
0 & \ldots \\ q^3 & q^3 & q^3 & q^3 & q^3 & 0 & \ldots \\ q^4 & q^4 & q^4
& q^4 & q^4 & q^4 & \ldots \\q^5 & q^5 & q^5 & q^5 & q^5 & q^5
&\ldots\\
\vdots &
\vdots & \vdots & \vdots & \vdots & \vdots &
\ddots\end{array}\right).  \end{displaymath} Let $\tilde{N}$ be the matrix generated by this production matrix. Dividing column $k$ of the $\tilde{M}$ by
$$\prod_{i=1}^k a_i =\prod_{i=1} q^i=q^{\binom{k}{2}},$$ we recover the matrix $N$.

The generating function $$\cfrac{1}{1-\cfrac{ x}{1-\cfrac{qx}{1-\cfrac{q^2 x}{1-\cdots}}}}$$ is equivalent to
$$\cfrac{1}{1-x-\cfrac{qx^2}{1-(q+q^2)x-\cfrac{q^5 x^2}{1-(q^3+q^4)x-\cfrac{q^9x^2}{1-(q^5+q^6)x-\cdots}}}}.$$
\noindent This leads to the production matrix
\begin{displaymath}\left(\begin{array}{ccccccc} 1 & 1 &
0
& 0 & 0 & 0 & \ldots \\q & q+q^2 & 1 & 0 & 0 & 0 & \ldots \\ 0 & q^5
& q^3+q^4 & 1 & 0 &
0 & \ldots \\ 0 & 0 & q^9 & q^5+q^6 & 1 & 0 & \ldots \\ 0 & 0 & 0
& q^{13} & q^7+q^8 & 1 & \ldots \\0 & 0 & 0 & 0 & q^{17} & q^9+q^{10}
&\ldots\\
\vdots &
\vdots & \vdots & \vdots & \vdots & \vdots &
\ddots\end{array}\right), \end{displaymath} which generates the matrix $M$, with first column equal to $\mu_n$.
\noindent The inverse of the matrix $M$ is the coefficient array of the orthogonal polynomials defined by
$$P_n(x)=(x-q^{2n-3}(1+q))P_{n-1}(x)-q^{4n-7}P_{n-2}(x),$$ where $P_0(x)=1$ and $P_1(x)=x-1$.

For $N^{-1} \cdot M$, we obtain the matrix that begins

\begin{displaymath}\left(\begin{array}{ccccccc} 1 & 0 &
0
& 0 & 0 & 0& \ldots \\0 & 1 & 0 & 0 & 0 & 0& \ldots \\ 0 & q^2
& 1 & 0 & 0 &
 0& \ldots \\ 0 & q^5 & q^3+q^4 & 1 & 0 & 0& \ldots \\ 0 & q^9 & q^7+q^8+q^9
& q^4+q^5+q^6& 1 & 0&  \ldots \\
0 & q^{14} & q^{12}+q^{13}+q^{14}+q^{15}
& q^9+q^{10}+q^{11}+q^{12}+q^{13}& q^5+q^6+q^7+q^8 & 1&  \ldots \\
\vdots &
\vdots & \vdots & \vdots & \vdots  & \vdots &
\ddots\end{array}\right). \end{displaymath}
\noindent Dividing each element $M_{n,k}$ of $M$ by
$$q^{\binom{n-k+2}{2}-1},$$ we obtain the matrix
\begin{displaymath}\left(\begin{array}{ccccccc} 1 & 0 &
0
& 0 & 0 & 0& \ldots \\0 & 1 & 0 & 0 & 0 & 0& \ldots \\ 0 & 1
& 1 & 0 & 0 &
 0& \ldots \\ 0 & 1 & q(1+q) & 1 & 0 & 0& \ldots \\ 0 & 1 & q^2(1+q+q^2)
& q^2(1+q+q^2) & 1 & 0&  \ldots \\
0 & 1 & q^3(1+q+q^2+q^3)
& q^4(1+q+2q^2+q^3+q^4) & q^3(1+q+q^2+q^3) & 1&  \ldots \\
\vdots &
\vdots & \vdots & \vdots & \vdots  & \vdots &
\ddots\end{array}\right), \end{displaymath} which is the Hadamard product of the matrices
\begin{displaymath}\left(\begin{array}{ccccccc} 1 & 0 &
0
& 0 & 0 & 0& \ldots \\0 & 1 & 0 & 0 & 0 & 0& \ldots \\ 0 & 1
& 1 & 0 & 0 &
 0& \ldots \\ 0 & 1 & 1+q & 1 & 0 & 0& \ldots \\ 0 & 1 & 1+q+q^2
& 1+q+q^2 & 1 & 0&  \ldots \\
0 & 1 & 1+q+q^2+q^3
& 1+q+2q^2+q^3+q^4 & 1+q+q^2+q^3 & 1&  \ldots \\
\vdots &
\vdots & \vdots & \vdots & \vdots  & \vdots &
\ddots\end{array}\right), \end{displaymath} and
\begin{displaymath}\left(\begin{array}{ccccccc} 1 & 0 &
0
& 0 & 0 & 0& \ldots \\0 & 1 & 0 & 0 & 0 & 0& \ldots \\ 0 & 1
& 1 & 0 & 0 &
 0& \ldots \\ 0 & 1 & q & 1 & 0 & 0& \ldots \\ 0 & 1 & q^2
& q^2 & 1 & 0&  \ldots \\
0 & 1 & q^3
& q^4 & q^3 & 1&  \ldots \\
\vdots &
\vdots & \vdots & \vdots & \vdots  & \vdots &
\ddots\end{array}\right), \end{displaymath} where the first matrix is the
$q$-Riordan array ${n-1 \brack n-k}_q$ \cite{qAnalogue}, and the second matrix is a shifted version of the matrix
$q^{k(n-k)}$.

The production matrix of $N^{-1}\cdot M$ begins
\begin{displaymath}\left(\begin{array}{ccccccc}
0 & 1 & 0 & 0 & 0 & 0 & \ldots \\
0 & q^2 & 1 & 0 & 0 & 0 & \ldots \\
0 & q^5(q-1) & q^2(q^2+q-1) & 1 & 0 &
0 & \ldots \\ 0 & 0 & q^7(q^2-1) & q^3(q^3+q^2-1) & 1 & 0 & \ldots \\ 0 & 0 & 0
& q^{10}(q^3-1) & q^4(q^4+q^3-1) & 1 & \ldots \\0 & 0 & 0 & 0 & q^{13}(q^4-1) & q^5(q^5+q^4-1)
&\ldots\\
\vdots &
\vdots & \vdots & \vdots & \vdots & \vdots &
\ddots\end{array}\right), \end{displaymath} indicating that in this case, the inverse matrix $(N^{-1}\cdot M)^{-1}=M^{-1}\cdot N$ is the coefficient array of a family of orthogonal polynomials whose parameters are given by the production matrix above.

We look more closely at the case of $q=2$. We find that
\begin{displaymath}M=\left(\begin{array}{ccccccc} 1 & 0 &
0
& 0 & 0 & 0& \ldots \\1 & 1 & 0 & 0 & 0 & 0& \ldots \\ 3 & 7
& 1 & 0 & 0 &
 0& \ldots \\ 17 & 77 & 31 & 1 & 0 & 0& \ldots \\ 171 & 1471 & 1333
& 127 & 1 & 0&  \ldots \\
3113 & 51653 & 98487
& 21717 & 511 & 1&  \ldots \\
\vdots &
\vdots & \vdots & \vdots & \vdots  & \vdots &
\ddots\end{array}\right), \end{displaymath} while
\begin{displaymath}N=\left(\begin{array}{ccccccc} 1 & 0 &
0
& 0 & 0 & 0& \ldots \\1 & 1 & 0 & 0 & 0 & 0& \ldots \\ 3 & 3
& 1 & 0 & 0 &
 0& \ldots \\ 17 & 17 & 7 & 1 & 0 & 0& \ldots \\ 171 & 171 & 77
& 51 & 1 & 0&  \ldots \\
3113 & 3113 & 1471
& 325 & 31 & 1&  \ldots \\
\vdots &
\vdots & \vdots & \vdots & \vdots  & \vdots &
\ddots\end{array}\right). \end{displaymath}
\noindent Then
\begin{displaymath}N^{-1}\cdot M=\left(\begin{array}{ccccccc} 1 & 0 &
0
& 0 & 0 & 0& \ldots \\0 & 1 & 0 & 0 & 0 & 0& \ldots \\ 0 & 4
& 1 & 0 & 0 &
 0& \ldots \\ 0 & 32 & 24 & 1 & 0 & 0& \ldots \\ 0 & 512 & 896
& 112 & 1 & 0&  \ldots \\
0 & 16384 & 61440
& 17920 & 480 & 1&  \ldots \\
\vdots &
\vdots & \vdots & \vdots & \vdots  & \vdots &
\ddots\end{array}\right). \end{displaymath}
\noindent Looking at the reduced matrix
\begin{displaymath}\left(\begin{array}{ccccccc} 1 & 0 &
0
& 0 & 0 & 0& \ldots  \\  4
& 1 & 0 & 0 &
 0& 0& \ldots \\ 32 & 24 & 1 & 0 & 0& 0&\ldots \\ 512 & 896
& 112 & 1 & 0& 0& \ldots \\
16384 & 61440
& 17920 & 480 & 1& 0& \ldots \\
1048576 & 8126464
& 5079040 & 317440 & 1984& 1& \ldots \\
\vdots &
\vdots & \vdots & \vdots & \vdots  & \vdots &
\ddots\end{array}\right) \end{displaymath} we see that it is the moment array of the family of orthogonal polynomials whose parameters are given in the production matrix
\begin{displaymath}\left(\begin{array}{ccccccc} 4 & 1 &
0
& 0 & 0 & 0& \ldots \\16 & 20 & 1 & 0 & 0 & 0& \ldots \\ 0 & 384
& 88 & 1 & 0 &
 0& \ldots \\ 0 & 0 & 7168 & 368 & 1 & 0& \ldots \\ 0 & 0 & 0
& 122880 & 1504 & 1 &  \ldots \\
0 & 0 & 0
& 0 & 2031616 & 6080&  \ldots \\
\vdots &
\vdots & \vdots & \vdots & \vdots  & \vdots &
\ddots\end{array}\right). \end{displaymath}  \noindent We deduce that the sequence
$1,4,32,512,16384,\ldots$ or
$2^{n(n+3)/2}$ \seqnum{A036442} has a generating function given by
$$\cfrac{1}{1-4x-\cfrac{16x^2}{1-20x-\cfrac{384x^2}{1-88x-\cfrac{7168x^2}{1-368x-\cdots}}}},$$ or equivalently,
$$\cfrac{1}{1-\cfrac{4x}{1-\cfrac{4x}{1-\cfrac{16x}{1-\cfrac{24x}{1-\cfrac{64x}{1-\cdots}}}}}}.$$
\noindent In this latter expression, the coefficients are given by the sequence
$$b(n)=2^{n+2}-2^{(n+1)/2}(1-(-1)^n).$$
The Hankel transform of $2^{n(n+3)/2}$ is then given by \cite{Kratt1, Kratt2, Layman}
$$h_n=\prod_{k=0}^{n-1} (b(2k+1)b(2k+2))^{n-k}.$$
A similar analysis can be carried out for $q^{n(n+3)/2}$.

\end{example}

\begin{example} \textbf{The little Schr\"oder numbers}. In this example, we take a base sequence $a_n$ given by
$$1,1,2,1,2,1,2,1,2,1,2,1,2,1,2,\cdots.$$ The sequence with generating function
$$\cfrac{1}{1-\cfrac{a_1x}{1-\cfrac{a_2x}{1-\cfrac{a_3x}{1-\cdots}}}}=\cfrac{1}{1-\cfrac{x}{1-\cfrac{2x}{1-\cfrac{x}{1-\cdots}}}}$$ is the sequence of little Schr\"oder numbers \seqnum{A001003}
$$1, 1, 3, 11, 45, 197, 903,\ldots.$$ \noindent These numbers are also generated by
$$\cfrac{1}{1-x-\cfrac{2x^2}{1-3x-\cfrac{2x^2}{1-3x-\cfrac{2x^2}{1-3x-\cdots}}}}.$$ \noindent
We have
\begin{displaymath}\left(\begin{array}{ccccccc} 1 & 0 &
0
& 0 & 0 & 0& \ldots  \\  -1
& 1 & 0 & 0 &
 0& 0& \ldots \\ 0 & -2 & 1 & 0 & 0& 0&\ldots \\ 0 & 0
& -1 & 1 & 0& 0& \ldots \\
0 & 0
& 0 & -2 & 1& 0& \ldots \\
 & 0
& 0 & 0 & -1 & 1& \ldots \\
\vdots &
\vdots & \vdots & \vdots & \vdots  & \vdots &
\ddots\end{array}\right)^{-1}=\left(\begin{array}{ccccccc} 1 & 0 &
0
& 0 & 0 & 0& \ldots  \\  1
& 1 & 0 & 0 &
 0& 0& \ldots \\ 2 & 2 & 1 & 0 & 0& 0&\ldots \\ 2 & 2
& 1 & 1 & 0& 0& \ldots \\
4 & 4
& 2 & 2 & 1& 0& \ldots \\
4 & 4
& 2 & 2 & 1 & 1& \ldots \\
\vdots &
\vdots & \vdots & \vdots & \vdots  & \vdots &
\ddots\end{array}\right), \end{displaymath} so that the matrix $N$ in this case begins
\begin{displaymath}\left(\begin{array}{ccccccc} \color{blue}{1} & 0 &
0
& 0 & 0 & 0& \ldots  \\  \color{blue}{1}
& \color{red}{1} & 0 & 0 &
 0& 0& \ldots \\ \color{blue}{3} & \color{red}{3} & \color{blue}{1} & 0 & 0& 0&\ldots \\ \color{blue}{11} & \color{red}{11}
& \color{blue}{4} & \color{red}{1} & 0& 0& \ldots \\
\color{blue}{45} & \color{red}{45}
& \color{blue}{17} & \color{red}{6} & \color{blue}{1} & 0& \ldots \\
\color{blue}{197} & \color{red}{197}
& \color{blue}{76} & \color{red}{31} & \color{blue}{7} & \color{red}{1}& \ldots \\
\vdots &
\vdots & \vdots & \vdots & \vdots  & \vdots &
\ddots\end{array}\right) \end{displaymath} with production matrix
\begin{displaymath}\left(\begin{array}{ccccccc} \color{blue}{1} & \color{red}{1} &
0
& 0 & 0 & 0& \ldots \\\color{blue}{2} & \color{red}{2} & \color{blue}{1} & 0 & 0 & 0& \ldots \\
\color{blue}{2} & \color{red}{2} & \color{blue}{1} & \color{red}{1} & 0 &
 0& \ldots \\
 \color{blue}{4} & \color{red}{4} & \color{blue}{2} & \color{red}{2} & \color{blue}{1} & 0& \ldots \\
 \color{blue}{4} & \color{red}{4} & \color{blue}{2} & \color{red}{2} & \color{blue}{1} & \color{red}{1} &  \ldots \\
\color{blue}{8} & \color{red}{8} & \color{blue}{4}
& \color{red}{4} & \color{blue}{2} & \color{red}{2} &  \ldots \\
\vdots &
\vdots & \vdots & \vdots & \vdots  & \vdots &
\ddots\end{array}\right). \end{displaymath} \noindent
\noindent
For instance, we have
$$\color{blue}{17}=\color{blue}{1}.\color{red}{11}\color{black}{+}\color{blue}{1}.\color{blue}{4}\color{black}{+}\color{blue}{2}.\color{red}{1}\color{black}{+}\color{blue}{2}.0\color{black}{+}\cdots,$$
and
$$\color{red}{45}=\color{red}{1}.\color{blue}{11}\color{black}{+}\color{red}{2}.\color{red}{11}\color{black}{+}\color{red}{2}.\color{blue}{4}\color{black}{+}\color{red}{4}.\color{red}{1}\color{black}{+}\cdots.$$
\noindent The matrix $M$ is given by
\begin{displaymath}\left(\begin{array}{ccccccc} 1 & 0 &
0
& 0 & 0 & 0& \ldots  \\  1
& 1 & 0 & 0 &
 0& 0& \ldots \\ 3 & 4 & 1 & 0 & 0& 0&\ldots \\ 11 & 17
& 7 & 1 & 0& 0& \ldots \\
45 & 76
& 40 & 10 & 1& 0& \ldots \\
197 & 353
& 216 & 72 & 13 & 1& \ldots \\
\vdots &
\vdots & \vdots & \vdots & \vdots  & \vdots &
\ddots\end{array}\right). \end{displaymath} \noindent This is the Riordan array \seqnum{A172094}
$$\left(\frac{1+x-\sqrt{1-6x+x^2}}{4x}, \frac{1-3x-\sqrt{1-4x+x^2}}{4x}\right)=\left(\frac{1}{1+x}, \frac{x}{1+3x+2x^2}\right)^{-1},$$ with production matrix
\begin{displaymath}\left(\begin{array}{ccccccc} 1 & 1 &
0
& 0 & 0 & 0& \ldots  \\  2
&3 & 1 & 0 &
 0& 0& \ldots \\ 0 & 2 & 3 & 1 & 0& 0&\ldots \\ 0 & 0
& 2 & 3 & 1& 0& \ldots \\
0 & 0
& 0 & 2 & 3& 1& \ldots \\
0 & 0
& 0 & 0 & 2 & 3& \ldots \\
\vdots &
\vdots & \vdots & \vdots & \vdots  & \vdots &
\ddots\end{array}\right). \end{displaymath}
\noindent The matrix $N$ is a ``mixture" (in left to right interleaved fashion) \cite{Davenport} of this Riordan array and the related Riordan array
$$\left(1, \frac{1-3x-\sqrt{1-6x+x^2}}{4x}\right)=\left(1, \frac{x}{1+3x+2x^2}\right)^{-1},$$
 which has production matrix
\begin{displaymath}\left(\begin{array}{ccccccc} 0 & 1 &
0
& 0 & 0 & 0& \ldots  \\  0
&3 & 1 & 0 &
 0& 0& \ldots \\ 0 & 2 & 3 & 1 & 0& 0&\ldots \\ 0 & 0
& 2 & 3 & 1& 0& \ldots \\
0 & 0
& 0 & 2 & 3& 1& \ldots \\
0 & 0
& 0 & 0 & 2 & 3& \ldots \\
\vdots &
\vdots & \vdots & \vdots & \vdots  & \vdots &
\ddots\end{array}\right). \end{displaymath}
We have
\begin{displaymath}
N_{n,k}=\begin{cases}
[x^n] \frac{1+x-\sqrt{1-6x+x^2}}{4x}\left(\frac{1-3x-\sqrt{1-6x+x^2}}{4}\right)^{k/2},\text{if $k$ is even}\\
[x^n] x^k \left(\frac{1-3x-\sqrt{1-6x+x^2}}{4x^2}\right)^{(k+1)/2}, \text{if $k$ is odd}.
\end{cases} \end{displaymath}
We note that in like fashion, the matrix $N^{-1}$, which begins
\begin{displaymath}\left(\begin{array}{ccccccc} \color{blue}{1} & 0 &
0
& 0 & 0 & 0& \ldots  \\  \color{red}{-1}
& \color{red}{1} & 0 & 0 &
 0& 0& \ldots \\ 0 & \color{blue}{-3} & \color{blue}{1} & 0 & 0& 0&\ldots \\ 0 & \color{red}{1}
& \color{red}{-4} & \color{red}{1} & 0& 0& \ldots \\
0 & 0
& \color{blue}{7} & \color{blue}{-6} & \color{blue}{1}& 0& \ldots \\
0 & 0
& \color{red}{-1} & \color{red}{11} & \color{red}{-7} & \color{red}{1} & \ldots \\
\vdots &
\vdots & \vdots & \vdots & \vdots  & \vdots &
\ddots\end{array}\right), \end{displaymath} is a ``mixture" (in shifted alternate row fashion) of the two matrices
$$\left(\frac{1}{1+3x+2x^2},\frac{x}{1+3x+2x^2}\right)\quad \text{and}\quad \left(\frac{1}{1+x}, \frac{x}{1+3x+2x^2}\right).$$
\noindent For instance, the array $\left(\frac{1}{1+3x+2x^2},\frac{x}{1+3x+2x^2}\right)$ begins
\begin{displaymath}\left(\begin{array}{ccccccc} 1 & 0 &
0
& 0 & 0 & 0& \ldots  \\  -3
& 1 & 0 & 0 &
 0& 0& \ldots \\ 7 & -6 & 1 & 0 & 0 & 0&\ldots \\ -15 & 23
& -9 & 1 & 0& 0& \ldots \\
31 & -72
& 48 & -12 & 1& 0& \ldots \\
-63 & 201
& -198 & 82 & -15 & 1& \ldots \\
\vdots &
\vdots & \vdots & \vdots & \vdots  & \vdots &
\ddots\end{array}\right), \end{displaymath}
 while the array $\left(\frac{1}{1+x}, \frac{x}{1+3x+2x^2}\right)$ begins

\begin{displaymath}\left(\begin{array}{ccccccc} 1 & 0 &
0
& 0 & 0 & 0& \ldots  \\  -1
& 1 & 0 & 0 &
 0& 0& \ldots \\ 1 & -4 & 1 & 0 & 0& 0&\ldots \\ -1 & 11
& -7 & 1 & 0& 0& \ldots \\
1 & -26
& 30 & -10 & 1& 0& \ldots \\
-1 & 57
& -102 & 58 & -13 & 1& \ldots \\
\vdots &
\vdots & \vdots & \vdots & \vdots  & \vdots &
\ddots\end{array}\right). \end{displaymath}

\noindent We have
\begin{displaymath}N^{-1}\cdot M=\left(\begin{array}{ccccccc} 1 & 0 &
0
& 0 & 0 & 0& \ldots  \\  0
& 1 & 0 & 0 &
 0& 0& \ldots \\ 0 & 1 & 1 & 0 & 0& 0&\ldots \\ 0 & 2
& 3 & 1 & 0& 0& \ldots \\
0 & 2
& 5 & 4 & 1& 0& \ldots \\
0 & 4
& 12 & 13 & 6 & 1& \ldots \\
\vdots &
\vdots & \vdots & \vdots & \vdots  & \vdots &
\ddots\end{array}\right). \end{displaymath}
We find that the production matrix of the inverse of this matrix is given by
\begin{displaymath}\left(\begin{array}{ccccccc} 0 & 1 &
0
& 0 & 0 & 0& \ldots \\0 & -1 & 1 & 0 & 0 & 0& \ldots \\ 0 & 0
& -2 & 1 & 0 &
 0& \ldots \\ 0 & 0 & 0 & -1 & 1 & 0& \ldots \\ 0 & 0 & 0
& 0 & -2 & 1 &  \ldots \\
0 & 0 & 0
& 0 & 0 & -1 &  \ldots \\
\vdots &
\vdots & \vdots & \vdots & \vdots  & \vdots &
\ddots\end{array}\right). \end{displaymath} \noindent This is the beheading of the inverse of the matrix
\begin{displaymath}\left(\begin{array}{ccccccc} 1 & 0 &
0
& 0 & 0 & 0& \ldots  \\  0
& 1 & 0 & 0 &
 0& 0& \ldots \\ 0 & 1 & 1 & 0 & 0& 0&\ldots \\ 0 & 2
& 2 & 1 & 0& 0& \ldots \\
0 & 2
& 2 & 1 & 1& 0& \ldots \\
0 & 4
& 4 & 2 & 2 & 1& \ldots \\
\vdots &
\vdots & \vdots & \vdots & \vdots  & \vdots &
\ddots\end{array}\right). \end{displaymath}
\noindent The form of the production matrix of the inverse is reflected in the structure of $N^{-1}\cdot M$ as follows: the internal elements of each even row satisfy $$t_{i,j}=1.t_{1-1,j-1}+1.t_{i-1,j},$$ while for odd rows we have
$$t_{i,j}=1.t_{i-1,j-1}+2.t_{i-1,j}.$$
\noindent We remark that it is clear that the interleaved structure of $N$, based on two Riordan arrays, will be replicated in the case of any sequence $a_n$ of the form $1,1,r,1,r,1,r,1,\ldots$.

\end{example}
\section{Conclusion} Using the parameters of equivalent Stieltjes and Jacobi continued fractions, we have defined two matrices $N$ and $M$, and we have studied the product $N^{-1}M$ in three specific cases. In each case, some noteworthy results have emerged. We conclude that the matrix $N^{-1}M$ is worthy of further study.

\bigskip
\hrule
\bigskip
\noindent 2010 {\it Mathematics Subject Classification}: Primary
15B36; Secondary 11B83, 15A09, 30B70, 42C05.

\noindent \emph{Keywords:} Matrix, Stieltjes continued fraction, Jacobi continued fraction, orthogonal polynomials, production matrix, Riordan array, Hankel transform

\bigskip
\hrule
\bigskip
\noindent Concerns sequences
\seqnum{A000108},
\seqnum{A001003},
\seqnum{A033184},
\seqnum{A036442},
\seqnum{A039599},
\seqnum{A172094}.

\end{document}